\documentclass[11pt]{article}
\usepackage{graphicx}
\usepackage{amsmath}
\usepackage{amsfonts}
\usepackage{epsfig}
\usepackage{setspace}

\newcommand{\be}{\begin{equation}}
\newcommand{\ee}{\end{equation}}
\newcommand{\beqn}{\begin{eqnarray}}
\newcommand{\eeqn}{\end{eqnarray}}
\newcommand{\beqns}{\begin{eqnarray*}}
\newcommand{\eeqns}{\end{eqnarray*}}

\newcommand{\Cov}{\mbox{Cov}\ }

\newcommand{\EE}{\ensuremath{{\mathbb E}}}
\newcommand{\II}{\ensuremath{{\mathbb I}}}

\newcommand{\fr}[1]{(\ref{#1})}

\newcommand{\om}{\omega}

\newcommand{\Om}{\Omega}

\newcommand{\bom}{\mbox{\mathversion{bold}$\om$}}

\newtheorem{lemma}{Lemma}
\newtheorem{theorem}{Theorem}

\newtheorem{remark}{Remark}

\begin{document}

\title{\Large{\bf Anisotropic functional Fourier deconvolution with long-memory dependent errors: a minimax study }}

\author{
\large{ Rida Benhaddou}  \footnote{E-mail address: Benhaddo@ohio.edu}
  \\ \\
Department of Mathematics, Ohio University, Athens, OH 45701} 
\date{}

\doublespacing
\maketitle
\begin{abstract}
We  investigate minimax results for the anisotropic functional deconvolution model when observations are affected by the presence of long-memory. Under specific conditions about the covariance matrices of the errors, we follow a standard procedure to construct an adaptive wavelet-based estimator that attains asymptotically near-optimal convergence rates. These rates depend on the parameter associated with the weakest long-range dependence, and deteriorate as the intensity of long-memory increases. This behavior suggests that the estimator adjusts to the best case scenario and that the weakest LM dominates.

{\bf Keywords and phrases: Anisotropic functional deconvolution,  Besov space, long memory, minimax convergence rates}\\ 

{\bf AMS (2000) Subject Classification: 62G05, 62G20, 62G08 }
 \end{abstract} 

\section{Introduction.}
Consider the problem of estimating the two-dimensional function $f(. , .)\in L^2(U)$ based on observations 
\be \label{inobserv}
Y(t_i, x_l)=\int^1_0f(s, x_l)g(t_i-s, x_l)ds+ \sigma \xi_{il}, \ \ i=1, 2, \cdots, N, \ l=1, 2, \cdots, M.
\ee
where $t_i=\frac{i}{N}$, $i=1, 2, \cdots, N$, $x_l=\frac{l}{M}$, $l=1, 2, \cdots, M$, $U=[0, 1]*[0, 1]$ and $\xi_{il}$ are Gaussian random variables that are independent for different $l=1, 2, \cdots, M$, but dependent for different $i=1, 2, \cdots, N$. The function $f(. , .)$ is periodic, and the convolution kernel is  known to the experimenter. Model \fr{inobserv} is referred to as functional deconvolution model. This model is motivated by experiments in which one needs to recover a two-dimensional function using  observations  of its convolutions along profiles $x=x_l$, $l=1, 2, \cdots, M$. This situation occurs, for example, 
in seismic inversions. Let $\xi^l_N$ be zero mean vector with components $\xi_{il}$,  $i=1, 2, \cdots, N$, and let $\Sigma^l_N=\Cov(\xi^l_N)=\EE\left[\xi^l_N(\xi^l_N)^T\right]$ be its covariance matrix. Consider the following assumption regarding the errors $\xi_{il}$ and their covariance matrices $\Sigma^l_N$.\\
\noindent
{\bf Assumption A.1.  }\label{A1}For each covariance matrix $\Sigma^l_N$, $l=1, 2, \cdots, M$, there exist constants $c_1$ and $c_2$ ($0<c_1 \leq c_2 < \infty$), independent of $N$, such that 
\be \label{assum1}
c_1 N^{1-\alpha_l}\leq \lambda_{\min}\left(\Sigma^l_N\right)\leq \lambda_{\max}\left(\Sigma^l_N\right)\leq c_2 N^{1-\alpha_l}, \ \ 0<\alpha_l \leq 1.
\ee
where $\alpha_l\in (0, 1]$, is the long-memory parameter associated with vector $\xi^l_N$, and $\lambda_{\min}\left(\Sigma^l_N\right)$ and $\lambda_{\max}\left(\Sigma^l_N\right)$ are the smallest and the largest eigenvalues of $\Sigma^l_N$, respectively.\\
\noindent
Assumption $A.1$ is valid for example when $\xi^l_N$ are fractional Gaussian noises or fractional ARIMA, (e.g., see Benhaddou et al.~(2014), Section 2,  for more detail).  

Deconvolution model has been the subject of a great deal of papers since late 1980s, but the most significant contribution was that of Donoho~(1995) who was the first to devise a wavelet solution to the problem. Other attempts include,  Abramovich and  Silverman (1998), Walter  and  Shen (1999), Donoho  and  Raimondo (2004), Johnstone et al.~(2004),  among others.  In the case of functional deconvolution model with $f(t, x)\equiv f(t)$, Pensky and sapatinas~(2009, 2010, 2011) pioneered into the formulation and further development of the problem. 
 
Functional deconvolution model of type \fr{inobserv} with $\alpha_1= \alpha_2=\cdots= \alpha_M=1$, corresponds to the $i.i.d.$ case studied in Benhaddou et al.~(2013). They constructed an adaptive hard-thresholding wavelet estimator, and showed that it is asymptotically near-optimal within a logarithmic factor of $MN$ under the $L^2$-risk over a wide range of anisotropic Besov balls. Benhaddou~(2017a) extends this work to the case of $L^p$-risk, $1\leq p <\infty$, to show that the hard-thresholding wavelet estimator is also asymptotically near-optimal. In these attempts, it is assumed that errors  are independent and identically distributed Gaussian random variables. However, empirical evidence has shown that, even at large lags, the correlation structure in the errors can decay at a power-like rate, rather than an exponential rate.
This phenomenon is referred to as long-memory (LM) or long-range dependence (LRD).

Long-memory  has been investigated quite considerably in the regression estimation framework, and to some less extent in the standard Fourier deconvolution model. In this latter case, one can list a few relevant endeavors;  Wang~(1997), Wishart~(2013), Benhaddou et al.~(2014), Kulik et al.~(2015) and Benhaddou~(2016). Recently, Benhaddou~(2017b) investigated the Laplace deconvolution problem for LM data when the unknown response function is defined on the positive real half-line.

The objective of the paper is to extend the work of Benhaddou et al.~(2013) to the case when the errors are affected by long-memory phenomenon. Like in Benhaddou et al.~(2014), the form of LM is not specified, rather, the covariance matrices of the errors are assumed to satisfy certain conditions in terms of their largest and smallest eigenvalues. Following a standard procedure, we derive minimax lower  bounds for the $L^2$-risk when $f(t , x)$ belongs to an anisotropic Besov ball and  the blurring function $g(t , x)$ is regular smooth. In addition, we show that the wavelet hard-thresholding estimator is adaptive and asymptotically near-optimal over a great array of Besov balls. Furthermore, we demonstrate that the proposed estimator attains convergence rates that depend on the parameter associated with the weakest LM from amongst the $M$ profiles. This suggests that the proposed estimator is affected by LM only to the extent of its weakest intensity among the profiles, and adjusts to the best case scenario. Besides, the long-memory phenomenon has a detrimental effect on the convergence rates. That is, the stronger the long-memory is, the slower the convergence rates will be, compared to the rates in Benhaddou et al.~(2013). 
 \section{Estimation Algorithm.}
 In what follows, denote the complex conjugate of $a$ by $\bar{a}$. 
  Consider a bounded bandwidth periodized wavelet basis (e.g., Meyer-type), $\psi_{j_1, k_1}(t)$, and a finitely supported periodized $s_o$-regular wavelet basis (e.g., Daubechies-type), $\eta_{j_2, k_2}(x)$. Both functions form orthonormal bases on $[0, 1]$, and therefore, the function $f(, . ,)$ can be expanded into a wavelet series as 
\be
f(t, x)= \sum_{\omega \in \Omega}\beta_{\omega}\psi_{j_1, k_1}(t)\eta_{j_2, k_2}(x).
\ee
where $\omega=\{j_1, k_1, j_2, k_2\}$, and 
\be \label{eq:setinf}
\Om = \left\{ \bom= (j_1,k_1; j_2,k_2): m_{i0} \leq j_i \leq \infty,\    k_i = 0, \cdots, 2^{j_i-1}; i=1,2  \right\}.
\ee 
Applying Fourier transform to equation \fr{inobserv} yields
\be \label{finobserv}
\tilde{Y}_{m}(x_l)=\tilde{f}_{m}(x_l)\tilde{g}_{m}(x_l)+\frac{\sigma}{\sqrt{N}}\tilde{\xi}_{m, l}, \ \ l=1, 2, \cdots, M.
\ee
where $\tilde{Y}_{m}(x_l)$, $\tilde{f}_{m}(x_l)$, $\tilde{g}_{m}(x_l)$ and $\tilde{\xi}_{m, l}$ are the Fourier coefficients of functions $Y(t_i, x_l)$, $g(t_i, x_l)$, $f(t_i, x_l)$ and $\xi_{il}$, respectively.  Let $\psi_{j_1, k_1, m}$ be Fourier coefficients of $\psi_{j_1, k_1}(t)$, then by Plancherel formula and \fr{finobserv}, an estimator for the wavelet coefficients is given by 
\be\label{tibeom}
\tilde{\beta}_{\omega}= \sum_{m\in W_{j_1}}\overline{\psi_{j_1, k_1,m}}\frac{1}{M}\sum^{M}_{l=1}\frac{\tilde{Y}_{m}(x_l)}{\tilde{g}_{m}(x_l)}\eta_{j_2, k_2}(x_l).
\ee
where, for any $j_1 \geq m_{10}$, 
\be \label{omeg}
W_{j_1}=\left \{ m: \psi_{j_1, k_1,m}\neq 0 \right\} \subseteq 2\pi/3\left[ -2^{j_1+2}, -2^{j_1} \right] \cup \left[ 2^{j_1}, 2^{j_1+2} \right ],
\ee
since Meyer wavelets are band-limited (see, e.g., Johnstone et al.~(2004)). Define 
\be \label{eq:setOmJ}
\Om(J_1, J_2) = \left\{ \bom= (j_1,k_1; j_2,k_2): m_{i0} \leq j_i \leq J_i-1,\    k_i = 0, \cdots, 2^{j_i-1}; i=1,2  \right\}.
\ee 
Then, allow the hard thresholding estimator for $f(t, x)$
\be \label{ef-hat}
\widehat{f}_{MN}(t, x)= \sum_{\omega \in \Omega(J_1, J_2)}\tilde{\beta}_{\omega} \II \left(|\tilde{\beta}_{\omega}| > \lambda^{MN}_{j_1} \right)\psi_{j_1, k_1}(t)\eta_{j_2, k_2}(x).
\ee
It remains to determine the choices of  $J_1$, $J_2$ and $\lambda^{MN}_{j_1}$ in \fr{ef-hat}. For that, it is necessary to evaluate the variance of \fr{tibeom}. Next is a condition that the convolution kernel $g(t, x)$ satisfies.

\noindent
{\bf Assumption A.2.  }\label{A1} The functional Fourier coefficients $g_m(x)$ of the kernel $g(t, x)$, for some positive constants $\nu$, and $K_1$ and $K_2$, independent of $m$ and $x$, are such that 
\be  \label{blur}
 K_1|m|^{-2\nu} \leq |g_m(x)|^2 \leq K_2 |m|^{-2\nu}.
\ee
Next we introduce a lemma which gives insight into the choice of the thresholds  $\lambda^{MN}_{j_1}$ and the maximal resolution levels $J_1$, $J_2$. \begin{lemma} \label{lem:Var}
Let $\tilde{ \beta}_{\omega}$ be defined in \fr{tibeom}. Then, under the conditions \fr{assum1} and \fr{blur}, one has 
\beqn 
\EE \left|\tilde{ \beta}_{\omega}-\beta_{\omega}\right|^{2} &\asymp&  \frac{\sigma^2}{M}2^{j_12\nu}\left[\frac{1}{M}\sum^M_{l=1}N^{-\alpha_l}|\eta_{j_2, k_2}(x_l)|^2\right].\label{var}\\
\EE|\tilde{\beta}_{\omega}-\beta_{\omega}|^4 &\asymp&  \frac{\sigma^42^{2j_1(2\nu+1)+j_2}}{M^3N^2}. \label{fourmomb}
\eeqn
\end{lemma}
According to Lemma \ref{lem:Var}, choose the thresholds $\lambda^{MN}_{j_1}$  as  
\be \label{thresh}
\lambda^{MN}_{j_1}= \rho 2^{\nu j_1}\left[ \frac{\ln(MN^{\alpha^*})}{MN^{\alpha^*}}\right]^{1/2},
\ee
where 
\be
\alpha^*=\max\{\alpha_1, \alpha_2, \cdots, \alpha_M\},
\ee
 and $J_1$ and $J_2$ as 
\be \label{maxres121}
2^{J_1}=\left[MN^{\alpha^*}\right]^{\frac{1}{2\nu +1}}, \ \ \ 2^{J_2}=MN^{\alpha^*},
\ee
where $\rho$ is some positive constant independent of $M$ and $N$. 
\section{Convergence rates and asymptotic optimality.}
Denote 
\beqn  \label{eq10}
 s^{*}_i&=&s_i+1/2 - 1/p,\ \ \ i=1, 2.\\
 s'_i&=&s_i + 1/2 -1/p'.
\eeqn
where $p'=\min\{2, p\}$. Assume that the function $f(t, x)$  and its wavelet coefficients satisfy the following.\\
\noindent
{\bf Assumption A.3.  } The function $f(t, x)$ belongs to an anisotropic two-dimensional Besov space. In particular, if $s_o \geq s_2$, its wavelet coefficients $ \beta_{\omega}$ satisfy
\be  \label{eq11}
 B^{s_1, s_2}_{p, q}(A)=\left \{ f \in L^2(U): \left( \sum_{j_1, j_2} 2^{(j_1s_1^{*}+j_2s_2^*)q}\left (\sum_{k_1, k_2}| \beta_{j_1, k_1, j_2, k_2}|^{p}\right)^{q/{p}}\right )^{1/q} \leq A\right \}.
\ee
To construct minimax lower bounds for the $L^2$-risk, we define the $L^2$-risk over the set  $\Theta$ as 
\be  \label{eq12}
 R_{MN}(\Theta)=\inf_{\tilde{f}_{MN}} \sup _{f \in \Theta}\EE \| \tilde{f}_{MN}-f\|^2,
\ee
where $\|g\|$ is the $L^2$-norm of a function $g$ and the infimum is taken over all possible estimators $\tilde{f}_{MN}$ of $f$. 

Notice that the thresholds $\lambda^{MN}_{j_1}$ and $J_1$ and $J_2$ are independent of the parameters $s_1$, $s_2$, $p$, $q$ and $A$ of the Besov ball $B^{s_1, s_2}_{p, q}(A)$, and therefore estimator \fr{ef-hat} is adaptive with respect to those parameters. It remains to see how such estimator performs asymptotically, so we shall evaluate the lower and upper bounds for the $L^2$-risk next. 

The following theorem provides the minimax lower bounds for the $L^2$-risk of any estimator $\tilde{f}_{MN}$.
\begin{theorem} \label{th:lowerbds}
Let $\min\{s_1, s_2\} \geq \max\{\frac{1}{p}, \frac{1}{2} \}$ with $1 \leq p,q \leq \infty$, and $A > 0$. Then, under conditions \fr{assum1}, \fr{blur} and \fr{eq11}, as $M, N \rightarrow \infty$, 
 \be \label{lowerbds}
R_{MN} (B^{s_1, s_2}_{p, q}(A)) \geq C A^2\left\{ \begin{array}{ll} 
 \left[\frac{\sigma^2}{A^2MN^{\alpha^*}} \right]^{\frac{2s_2}{2s_2+1}} , & \mbox{if}\ \  s_1 > s_2(2\nu + 1),\\
 \left[ \frac{\sigma^2}{A^{2}MN^{\alpha^*}}\right]^{\frac{2s_1}{2s_1 +2\nu +1}}   , & \mbox{if}\ \ \  (2\nu + 1)(\frac{1}{p}-\frac{1}{2})\leq s_1\leq s_2(2\nu + 1),\\
  \left[\frac{\sigma^2}{A^2MN^{\alpha^*}} \right]^{\frac{s^*_1}{2s^*_1+2\nu}}, &  \mbox{if}\  s_1< (2\nu + 1)(\frac{1}{p}-\frac{1}{2}).
\end{array} \right.
\ee
  \end{theorem}
{\bf Proof of Theorem \ref{th:lowerbds}}. In order to prove the theorem, we consider two cases, the case
when $f(t,x)$ is dense in both variables (the dense-dense case) and the case when $f(t,x)$ is dense in $x$ and sparse in $t$ (the sparse-dense case). Lemma $A.1$ of Bunea et al.~(2007) is then applied to find such lower bounds using conditions \fr{assum1}, \fr{blur} and \fr{eq11}. To complete the proof, we choose the highest of the lower bounds. $\Box$ 

The next lemma provides large deviation results for the wavelet coefficients $\beta_{\omega}$ and their estimates $\tilde{ \beta}_{\omega}$.
\begin{lemma} \label{lem:Lar-D} 
Let $\tilde{\beta}_{\omega}$ and $\lambda^{MN}_{j_1}$ be defined in $(7)$ and $(14)$, respectively. Let conditions \fr{assum1} and \fr{blur} hold. Then, for some positive constant $\gamma$, as $M$, $N\rightarrow \infty$, one has 
\be
\Pr\left(|\tilde{\beta}_{\omega}-\beta_{\omega}| > \gamma \lambda^{MN}_{j_1}\right)=O\left(\frac{\left[\frac{1}{MN}\right]^{\frac{\gamma^2 \rho^2}{2\sigma^2_o}}}{\sqrt{\ln(MN)}}\right)
\ee
where $\rho$ is defined in \fr{thresh} and $\sigma_o^2=\frac{c_2\sigma^2}{K_1}\left(\frac{8\pi}{3}\right)^{2\nu}$.
\end{lemma}
The next theorem gives the upper bounds for the minimax risk of the estimator \fr{ef-hat}.
\begin{theorem} \label{th:upperbds}
Let $\widehat{f}(. , .)$ be the wavelet estimator in \fr{ef-hat}, with $\lambda^{MN}_{j_1}$ given by \fr{thresh} and, $J_1$ and $J_2$ given by \fr{maxres121}. Let $s_0 > s_2$, and $\min\{s_1, s_2\} \geq \max\{\frac{1}{p}, \frac{1}{2} \}$, and let conditions \fr{assum1},  \fr{blur} and \fr{eq11} hold. If $\rho$ in \fr{thresh} is large enough, then, as $M, N \rightarrow \infty$, 
 \be \label{upperbds}
R_{MN}( B^{s_1, s_2}_{p, q}(A)) \leq C A^2\left\{ \begin{array}{ll} 
  \left[\frac{\sigma^2}{A^2MN^{\alpha^*}} \right]^{\frac{2s_2}{2s_2+1}}\left[\ln(MN)\right]^{\xi_1+\frac{2s_2}{2s_2+1}} , & \mbox{if}\  \ s_1 \geq s_2(2\nu + 1),\\
   \left[\frac{ \sigma^{2}}{A^2MN^{\alpha^*}} \right]^{\frac{2s_1}{2s_1 +2\nu +1}}\left[\ln(MN)\right]^{\frac{2s_1}{2s_1 +2\nu +1}} , & \mbox{if} \ \  (2\nu + 1)(\frac{1}{p}-\frac{1}{2})< s_1< s_2(2\nu + 1),\\
 \left[\frac{ \sigma^{2}}{A^2MN^{\alpha^*}} \right]^{\frac{2s'_1}{2s'_1+2\nu}}\left[\ln(MN)\right]^{\xi_2+{\frac{2s'_1}{2s'_1+2\nu}}}, &  \mbox{if}\ \  s_1\leq (2\nu + 1)(\frac{1}{p}-\frac{1}{2}).
\end{array} \right.
\ee
where $\xi_1$ and $\xi_2$ are defined as 
\be \label{d}
\xi_1 =\II \left( s_1= s_2(2\nu +1) \right), \ \ \
\xi_2 =\II  \left( s_1= (2\nu +1)\left(\frac{1}{p} -\frac{1}{2}\right) \right).
\ee
\end{theorem}
{\bf The proof of Theorem \ref{th:upperbds}.} The proof is very similar to that of Theorem 2 in Benhaddou et al.~(2013). $\Box$
\begin{remark}	
{\rm{
{\bf{(i)}}\,
 Theorems \ref{th:lowerbds} and \ref{th:upperbds} imply that, for the $L^2$-risk, the estimator \fr{ef-hat} is asymptotically near-optimal within a logarithmic factor of $MN$, over a wide range of anisotropic Besov balls $ B^{s_1, s_2}_{p, q}(A)$. \\
  \noindent {\bf{(ii)}}\,
 Notice that the rates of convergence are expressed in terms of the largest LM parameter $\alpha^*=\max\{\alpha_1, \alpha_2, \cdots, \alpha_M\}$, from amongst the $M$ profiles, which corresponds to the weakest long-range dependence. This implies that our estimator adjusts to the best case scenario and that the weakest LM dominates. \\
  \noindent {\bf{(iii)}}\, 
The convergence rates deteriorate as long-memory phenomenon gets stronger. More specifically, the stronger the LM is, the slower the convergence rates will be, compared to Benhaddou et al.~(2013).  This detrimental effect of LM on convergence rates was pointed out in Wishart~(2013), Kulik et al.~(2015), Benhaddou~(2016) and Benhaddou~(2017b).\\
 \noindent {\bf{(iv)}}\,
  For $\alpha_1=\alpha_2= \cdots=\alpha_M=1$, our rates of convergence match exactly those in Benhaddou et al.~(2013), by setting $\varepsilon^2= \frac{\sigma^2}{MN}$, and with $p=2$, those in Benhaddou~(2017a). \\
 \noindent {\bf{(v)}}\,
  Note that our rates of convergence are not directly comparable to those in Kulik et al.~(2015), since their rates pertain to an estimator of a one-dimensional function using a finite number of ($M$ different) channels, while in the present work the convergence rates are associated with the estimation of a two-dimensional function and that the number of profiles $M$ is asymptotic. 
    }}
\end{remark}
\section{Proofs. }
\subsection{Proof of the lower bounds.}
\underline{\bf  The dense-dense case.}
Using the same test functions $f_{\tilde{\omega}}$ and $f_{\omega}$, as in Benhaddou et al.~(2013), it can be shown that the $L^2$ norm of the difference satisfies    
\be  \label{Norm}
\|f_{\tilde{\omega}}-f_{\omega}\|^2_2 \geq \rho^2_{j_{1}j_2}  2^{j_1+j_2}/8,
\ee
In order to apply Lemma A.1 of Bunea et al.~(2007), one needs to verify condition $(ii)$. Denote $Q^{(N)}_{l, \omega}$ and $Q^{(N)}_{l, \tilde{\omega}}$, the vectors with components
\beqn
q_{\omega}(t_i, x_l)&=&g(t_i-., x_l)*f_{\omega}(., x_l), \ \ \ i=1, 2, \cdots, N.\\
q_{\tilde{\omega}}(t_i, x_l)&=&g(t_i-., x_l)*f_{\tilde{\omega}}(., x_l), \ \ \ i=1, 2, \cdots, N.
\eeqn
Then, the Kullback divergence is 
\beqn
K(P_{f_{\omega}}, P_{f_{\tilde{\omega}}})&=&\frac{1}{2\sigma^2}\sum^M_{l=1}\left(Q^{(N)}_{l, \omega}- Q^{(N)}_{l, \tilde{\omega}}\right)^T(\Sigma^l_N)^{-1} \left(Q^{(N)}_{l, \omega}- Q^{(N)}_{l, \tilde{\omega}}\right)\nonumber\\
&\leq& \frac{1}{2\sigma^2}\sum^M_{l=1}\lambda_{\max}\left[\left(\Sigma^l_N\right)^{-1}\right]\|Q^{(N)}_{l, \omega}- Q^{(N)}_{l, \tilde{\omega}}\|^2\nonumber\\
&\leq& \frac{N\rho^2_{j_1, j_2}2^{j_1+j_2}}{2\sigma^2c_1}\sum_{m\in W_{j_1}}\sum^M_{l=1}\lambda_{\max}\left(\Sigma^l_N\right)^{-1}|\tilde{g}_m(x_l)|^2|\psi_{j_1, k_1, m}|^2|\eta_{j_2, k_2}(x_l)|^2\nonumber\\
&\leq& \frac{MN^{\alpha^*}K_2}{2\sigma^2c_1}2 \pi\rho^2_{j_1, j_2}2^{j_1+j_2}2^{-2\nu j_1}\frac{1}{M}\sum^M_{l=1}|\eta_{j_2, k_2}(x_l)|^2,
\eeqn
Lemma A.1 suggests to choose $j_1$ and $j_2$ such that 
\be
 \frac{MN^{\alpha^*}K_2}{2\sigma^2c_1}2 \pi\rho^2_{j_1, j_2}2^{j_1+j_2}2^{-2\nu j_1} \leq 2^{j_1+j_2}\frac{\ln(2)}{16},
 \ee
Hence, using argument similar to Benhaddou et al.~(2013), the lower bounds are 
\be  \label{delta1}
\delta^2 =C A^2 \left\{ \begin{array}{ll}
\left[\frac{\sigma^{2}}{A^{2}MN^{\alpha^*}}\right]^{\frac{2s_1}{2s_1+2\nu +1}}, & \mbox{if}\ \ s_1 \leq s_2 (2\nu +1),\\
\left[\frac{\sigma^{2}}{A^{2}MN^{\alpha^*}}\right]^{\frac{2s_2}{2s_2 +1}},&  \mbox{if}\ \  s_1 > s_2 (2\nu +1).
\end{array} \right.
\ee
\underline{\bf  The sparse-dense case.}
Using the same test functions $f_{\tilde{\omega}}$ and $f_{\omega}$, as in Benhaddou et al.~(2013), and following the same procedure as in the dense-dense case,  it can be shown that the lower bounds are 
\be \label{delta2}
\delta^2 = C A^2\left\{ \begin{array}{ll}
\left[\frac{\sigma^{2}}{A^{2}MN^{\alpha^*}}\right]^{\frac{2s_2}{2s_2 +1}},&  \mbox{if}\ \  s^*_1 \geq s_2 2\nu ,\\
\left[\frac{\sigma^{2}}{A^{2}MN^{\alpha^*}}\right]^{\frac{2s^*_1}{2s^*_1+2\nu }},&  \mbox{if}\ \  \ s^*_1 < s_2 2\nu.
\end{array} \right.
\ee
To complete the proof, notice that the highest of the lower bounds corresponds to 
\be
d=\min \left\{ \frac{2s_1}{2s_1 + 2\nu +1}, \frac{2s_2}{2s_2 +1}, \frac{2s^*_1}{2s^*_1+2\nu } \right\}.
\ee
$\Box$
\subsection{Proof of the upper bounds.}
{\bf Proof of Lemma \ref{lem:Var}.} Note that 
\be \label{norrv}
\tilde{\beta}_{\omega}-\beta_{\omega}=\frac{\sigma}{\sqrt{N}}\sum_{m\in W_{j_1}}\overline{\psi_{j_1, k_1, m}}\frac{1}{M}\sum^M_{l=1}\frac{\tilde{\xi}_{m, l}}{\tilde{g}_m(x_l)}\eta_{j_2, k_2}(x_l),
\ee
Define the vector $U_{l}$,  with components $U_{m, l}=\overline{\psi_{j_1, k_1, m}}\frac{\eta_{j_2, k_2}(x_l)}{\tilde{g}_m(x_l)}$. Then,
\be \label{unorm}
\|U_{l}\|^2=\sum_{m\in W_{j_1}}|\psi_{j_1, k_1, m}|^2\frac{|\eta_{j_2, k_2}(x_l)|^2}{|\tilde{g}_m(x_l)|^2},
\ee
Hence, by conditions \fr{assum1} and \fr{blur}, and the fact that $|\psi_{j_1, k_1, m}|\leq 2^{-j_1/2}$, the variance of \fr{norrv} becomes
\beqn
\EE|\tilde{\beta}_{\omega}-\beta_{\omega}|^2&=&\frac{\sigma^2}{M^2N}\sum^M_{l=1}U^T_l\left(\Sigma^{(l)}_N\right)U_l\nonumber\\
&\leq& \frac{\sigma^2}{M^2N}\sum^M_{l=1}\lambda_{\max}\left(\Sigma^{(l)}_N\right)\|U_l\|^2\nonumber\\
&\leq& \frac{c_2\sigma^2}{M^2N}\sum^M_{l=1}N^{1-\alpha_l}\sum_{m\in W_{j_1}}|\psi_{j_1, k_1, m}|^2\frac{|\eta_{j_2, k_2}(x_l)|^2}{|\tilde{g}_m(x_l)|^2}\nonumber\\
&\leq& \frac{c_2\sigma^2}{K_1M}2^{2\nu j_1}\left[\frac{1}{M}\sum^M_{l=1}N^{-\alpha_l}|\eta_{j_2, k_2}(x_l)|^2\right],\label{35}
\eeqn
this completes the proof of \fr {var}. $\Box$

To prove \fr{fourmomb}, notice that the fourth moment of \fr{norrv} can be written as 
\beqn
\EE|\tilde{\beta}_{\omega}-\beta_{\omega}|^4&=&O\left( \frac{\sigma^4}{N^2M^4}\sum^M_{l=1}\left[\sum_{m\in W_{j_1}}\frac{|\psi_{j_1, k_1, m}|}{|\tilde{g}_m(x_l)|}\left(\EE|\tilde{\xi}_{m, l}|^4\right)^{\frac{1}{4}}\right]^4|\eta_{j_2, k_2}(x_l)|^4\right)\nonumber\\
&+& O\left(\left[\EE|\tilde{\beta}_{\omega}-\beta_{\omega}|^2\right]^2\right)\nonumber\\
&=& O\left( \frac{\sigma^4}{N^2M^4}\sum^M_{l=1}\left[\sum_{m\in W_{j_1}}\frac{|\psi_{j_1, k_1, m}|^2}{|\tilde{g}_m(x_l)|^2}\sum_{m\in W_{j_1}}\left(\EE|\tilde{\xi}_{m, l}|^2\right)\right]^2|\eta_{j_2, k_2}(x_l)|^4\right)\nonumber\\
&+& O\left(\left[\EE|\tilde{\beta}_{\omega}-\beta_{\omega}|^2\right]^2\right),
\eeqn
Now, since $\sum_{m\in W_{j_1}}\EE|\tilde{\xi}_{m, l}|^2 \asymp 2^{j_1}$, using $|\psi_{j_1, k_1, m}|\leq 2^{-j_1/2}$, conditions $(2)$ and $(4)$, completes the proof.  $\Box$\\
{\bf Proof of Lemma \ref {lem:Lar-D}.} Notice that the quantities $\tilde{\theta}_{\omega}=\tilde{\beta}_{\omega}-\beta_{\omega}$ are centered Gaussian random variables with variances of order \fr{35}. Hence, applying the Gaussian tail probability inequality completes the proof. $\Box$

\end{document}